\documentclass[12pt,a4paper]{article}
\usepackage{amsmath,amssymb,amsthm}
\usepackage{graphicx}
\usepackage{color}

\usepackage[notref,notcite]{showkeys}

\theoremstyle{plain}

\newtheorem{lemma}{Lemma}[section]

\theoremstyle{remark}
\newtheorem{remark}[lemma]{Remark}

\theoremstyle{definition}

\begin{document}

\newcommand{\strain}{{\boldsymbol\varepsilon}}
\newcommand{\straind}{{\boldsymbol e}}
\newcommand{\straine}{{\boldsymbol\varepsilon^e}}
\newcommand{\plast}{{\boldsymbol\varepsilon^p}}
\newcommand{\stress}{{\boldsymbol\sigma}}
\newcommand{\inter}{{\boldsymbol\xi}}
\newcommand{\ttau}{{\boldsymbol h_d}}
\newcommand{\et}{{\boldsymbol d}^{tr}}
\newcommand{\uu}{{\boldsymbol u}}                           
\newcommand{\no}{\nonumber}
\newcommand{\elena}{\color{red}}

                                %
\newcommand{\diffe}{D}
\title{\Large A new 3D macroscopic model for shape memory alloys describing
martensite reorientation}

                                %

                                %
\author{\normalsize Ferdinando Auricchio \\
{\small \it Dipartimento di Meccanica Strutturale} \\
{\small \it Universit\`a di Pavia, via Ferrata 1, 27100 Pavia, Italy}\\
{\scriptsize e-mail address: \tt auricchio@unipv.it} \\
\normalsize  Elena Bonetti\\
{\small \it Dipartimento di Matematica} \\
{\small \it Universit\`a di Pavia, via Ferrata 1, 27100 Pavia, Italy}\\
{\scriptsize e-mail address: \tt elena.bonetti@unipv.it}  }
\vspace{-1.5mm}

\date{\today}     
                                %
                                %
\maketitle
                                %
\begin{abstract}
In this paper we introduce a 3D phenomenological model for shape memory behavior, accounting for: martensite 
reorientation, asymmetric response of the material to tension/compression, different kinetics between forward and reverse phase transformation. 
We combine two modeling approaches using scalar and tensorial internal variables. Indeed,
we use volume proportions of different configurations of the crystal lattice 
(austenite and two variants of martensite) as scalar internal variables and 
the preferred direction of stress-induced martensite as tensorial internal variable. 
Then, we derive evolution equations by a generalization of the principle of virtual powers, including 
microforces and micromovements responsible for  phase transformation. 
In addition, we prescribe an evolution law for phase proportions  ensuring different evolution laws
during forward and reverse transformation of the oriented martensite.

\end{abstract}
\vskip3mm

   \renewcommand{\theequation}{\thesection.\arabic{equation}}                             
\section{Introduction}

In the last years shape memory alloys have been deeply investigated, 
from the point of view of modeling, analysis, and computation.
Indeed, these materials present many important industrial applications 
(e.g. aeronautical, biomedical, structural, and earthquake engineering) due to their characteristic 
of superelasticity
and shape memory effect.  

It is known that the shape memory effect is the consequence of a (reversible) martensitic 
phase transformation between different configurations of the crystal lattice in the alloy: 
from a high symmetric phase, austenite, to a lower symmetric configuration, martensite. 
Austenite is a solid phase (present at high temperatue) which can transform in martensite 
by means of a shearing mechanism. When transformation comes from thermal actions (lowering temperature) 
the result is a multi-direction martensite, in which variants compensate each other and there is no 
resulting macroscopic deformation. On the contrary, when transformation is obtained by loading, 
oriented martensite is formed in the stress direction, exhibiting a macroscopic deformation.

In particular, the research have been developed towards
the aim of finding a flexible phenomenological model. Some reliable models have been proposed 
to predict the response of such materials. 
Among the others, we focus on two models developed in the framework of phase transitions. 

The first, proposed by Fr\'emond (cf., e.g., \cite{fremondlibro} and some generalizations 
\cite{bonetti1}, \cite{bfl}) describes the behavior of shape memory in terms of (local) 
volume proportions of different configurations of the crystal lattice. More precisely, 
the austenite and two variants of martensite  are taken into  account. Note that 
the average behavior of different configurations 
is considered as the behavior of the equivalent single variant. The resulting model is able 
to describe phase transformation between different configurations. However, the model is 
obtained assuming that the direction of the transformation strain (associated to the detwinned martensite) is known. 

The second model we are considering has been proposed in \cite{souza} and then generalized in  
\cite{auricchiopet}. In this case one internal variable describes the phenomenon: the transformation strain tensor. 
In this case, the internal variable describes the direction of martensite orientation.  In addition, it leads to a simple and robust algorithm, based on plasticity-like return map.  
Thanks to this property the model has been used for implementation within finite element codes, allowing the simulation 
of complex SMA devices. However, some secondary effects are not included in this second model, as scalar and directional informations are tightly interconnected. 

Thus, one could wonder 
how to get a deeper description of micro-phenomena, possibly combining the main features of the two different 
approaches. 
Accordingly,
the purpose of this paper is to  answer to this idea combining 
the two theories describing secondary effects in the phase transitions 
as well as directional information for the transformation strain.
Thus, both scalar and tensorial internal variables are introduced, 
accounting for the phase proportions (assuming that in each point 
the phases may coexist with different proportions) and for the orientation 
of the transformation strain associated to the detwinned martensite. 
We consider both proportion and direction 
as internal variables and we write evolution equations for both of them. We recall that an attempt 
in this direction has been performed in  \cite{ARaltri}. 
In particular, we prescribe an evolution law to capture asymmetric 
response of the material in tension-compression loading. In \cite{ARS} 
a 1D model has been introduced to describe this kind of phenomenon, 
using an asymmetric energy  depending on a tensorial variable. 
However, this model seems to be hard to be extended to higher dimensions. On the contrary, our approach 
(developing an asymmetric evolution theory for phase proportions) can  apply to any space dimension.

\section{The model}

In this section we detail the derivation of the model. We mainly
refer to the approach proposed by Fr\'emond to describe the behavior of a  
thermomechanical system  in terms of state and dissipative variables, as well as energy and dissipation 
functionals (see \cite{fremondlibro}).
The main idea  consists in assuming that (microscopic) phase transformations are due to  micro-forces and micro-movements that have to be included 
in the global energy balance of the system (i.e. generalizing the principle of virtual powers). 
In particular, the equations governing 
the evolution of internal variables are recovered as balance equations (as for the momentum balance).

\subsection{The state and dissipative variables}

As it is known,  phase transformations in the alloy are due to the phase transitions occurring in 
the microstructure configuration between austenite  and twinned or detwinned martensite. 
In particular, detwinning manifests itself mainly through a shear strain so that we introduce a symmetric and deviatoric strain  which appears in presence of the detwinned martensite. The (local) volume proportions of austenite and martensite variants is represented by phase parameters
\begin{equation}\label{varfase}
\chi_A,\chi_M,\chi_S\in[0,1],\quad\chi_A+\chi_M+\chi_S=1.
\end{equation}
More precisely,
 $\chi_A$ stands for austenite, $\chi_M$ for twinned martensite, and $\chi_S$ for detwinned martensite. 
Furthermore, $\et$ is the direction of the deviatoric strain tensor associated to the detwinned martensite with $\|\et\|=\xi_s$
($\xi_s$ is the maximum amount for the detwinned martensite). 
Indeed, the deviatoric strain for detwinned martensite is given by $\chi_S\et$. 
Then,
$\theta$ is the absolute temperature,
$\strain(\uu)$ the (symmetric) linearized strain tensor ($\uu$ is the vector of small displacements as  
we restrict ourselves to small deformations).  Due to the internal constraint \eqref{varfase} on the phase proportions (coming by the 
their physical meaning) 
we can restrict ourselves to consider just two independent phase variables  $(\chi_M,\chi_S)$  letting 
$$
\chi_A=1-\chi_M-\chi_S,
$$
where 
$$
0\leq\chi_M,\chi_S\leq 1,\qquad \chi_M+\chi_S\leq 1.
$$
Finally, let us use the notation $\straine$ for the elastic component of the strain, so that it results
$$
\strain=\straine+\chi_S\et.
$$
Hence, the corresponding deviatoric strain ${\bf e}$ is
$$
{\bf e}:=\strain-\frac 1 3 \hbox{tr }(\strain){\bf I}
$$
${\bf I}$ being the identity matrix and $\hbox{tr }(\cdot)$ the trace operator. 
If $\stress$ is the Cauchy stress tensor, the deviatoric stress tensor ${\bf S}$ is
$$
{\bf S}:=\stress-\frac 1 3\hbox{tr }(\stress) {\bf I}=\stress-\sigma_m{\bf I}.
$$  
As far as evolution, this is described by dissipative variables $\chi_{Mt}$, $\chi_{St }$, $\et_t$, and $\nabla\theta$. These variables are in particular related to micro-velocities in the phase transformation. 

\begin{remark}
Let us comment about the choice of state variables. The main idea  consists in distinguishing between the norm and the direction of the inelastic strain.
In this way, we are able to describe the presence of a product phase,
 to which a homogenized strain is associated, and a parent phase, in which we find only elastic strain. However, in the parent phase, we can also distinguish between the presence of twinned martensite and austenite. 
Thus, we get a more complex and free description of the phenomenon with respect to
the Souza and Fr\'emond models, which could be useful in some situations.
\end{remark}

\subsection{The energy and dissipation functionals}

We introduce the following free energy functional (depending on state variables) as a combination of the energies associated to the single variants (combined with suitable proportions) and by an interaction energy, accounting also for internal constraints
\begin{equation}
\Psi(\strain,\et,\chi_M,\chi_S,\theta)=\Psi_{el}+\Psi_{id}+\Psi_{ch}+\Psi_v
\end{equation}
where
\begin{align}
&\Psi_{el}=\left(\frac \lambda 2+\frac \mu 3\right)(\hbox{tr }\strain)^2+\mu\|{\bf e}-\chi_S\et\|^2\\\no
&\Psi_{id}=c_s((\theta-\theta_0)-\theta\log\theta)\\\no
&\Psi_{ch}=(1-\chi_M-\chi_S)h_A(\theta)+\chi_M h_M(\theta)+\chi_S h_S(\theta)+\ttau(\theta):\et\\\no
&\Psi_v=I_K(\chi_M,\chi_S)+I_{\xi_s}(\|\et\|)+{\Psi_{int}}(\chi_M,\chi_S).
\end{align}
Here $c_s>0$ is the specific heat, $\lambda$ and $\mu$ are the Lam\'e constants; $h_A,h_S,h_M,{\bf h}_d$  are smooth thermal functions whose regularity  will be specified later on  (to ensure compatibility with thermodynamics). 
The function $I_K$ is the indicator function of the convex set $K$ 
$$
K:=\{(\chi_M,\chi_S)\in{\bf R}^2:0\leq\chi_M,\chi_S\leq1,\chi_M+\chi_S\leq 1\},
$$
i.e. it is
$I_K(\chi_M,\chi_S)=0$ if $(\chi_M,\chi_S)\in K$, while $I_K(\chi_M,\chi_S)=+\infty$ otherwise.
The function $I_{\xi_s}$ forces $\|\et\|=\xi_s$. Indeed, it is $I_{\xi_s}(\|\et\|)=0$ if $\|\et\|=\xi_s$ and it is $+\infty$ otherwise.  ${\Psi_{int}}$ is a (sufficiently) smooth function accounting for interaction energy. 
As  a possible choice for the interaction energy ${\Psi_{int}}$, we could simply consider 
\begin{eqnarray}
&\Psi_{int}(\chi_M,\chi_S)=C_{MS}\chi_M\chi_S+(C_{AM}\chi_M+C_{AS}\chi_S)(1-\chi_M-\chi_S)\\\no
&+C_{AMS}\chi_M\chi_S(1-\chi_M-\chi_S),
\end{eqnarray}
where $C_{MS},C_{AM},C_{AS}, C_{AMS}$ are positive constants.  

\begin{remark}
Note that for $(\chi_M,\chi_S)$ we have introduced a convex constraint forcing 
$(\chi_M,\chi_S)\in K$. The constraint on $\et$ is convex w.r.t. to its norm as it is $\|\et\|=\xi_s$. 
\end{remark}

Now, let us introduce the pseudo-potential of dissipation, which is a positive convex functional depending on dissipative variables, vanishing for vanishing dissipation (cf. \cite{moreau}). We have
\begin{equation}\label{pseudo}
 \phi(\chi_{Mt},\chi_{St},\et_t,\nabla\theta)=|\chi_{Mt}|+\phi_S(\chi_S,\stress,\chi_{St})
+\chi_S\|\et_t\|+\frac 1 {2\theta}|\nabla\theta|^2.
\end{equation}
  Note that, $\phi$ is considered 
to possibly ensure an asymmetric behavior in tension and compression and for forward and backward  transformation. 
This is due to the choice of the function $\phi_S$ (noting that it possibly depends on the stress and $\chi_S$). 
Indeed, this choice generalizes the classical situation for rate-independent systems, where it is
\begin{equation}
\phi_S(\chi_S,\stress,\chi_{St})=|\chi_{St}|.
\end{equation}
 Actually, $\phi_S$ is required to be rate independent with respect to 
$\chi_{St}$. Hence, accounting for a possible dependence in the evolution on the stress 
(e.g., for tension-compression behavior) and on the  volume of already detwinned martensite, we get
as  a further possible example
\begin{equation}\label{specifico}
\phi_S(\chi_S,\stress,\chi_{St})=d(\chi_S,\stress)(\chi_{St})^++|\chi_{St}|
\end{equation}
where 
$(f)^+=f$ if $f\geq 0$ and $(f)^+=0$ if $f\leq0$ and $d$ is a sufficiently smooth function. From now on we deal in particular with \eqref{specifico}.

\begin{remark}\label{asim}
 Note that we could refine the model, e.g. adding in \eqref{pseudo} a term as $\widehat d(\chi_S,\stress)(\chi_{St})^-$ for decreasing evolution of the product phase.
\end{remark}


\subsection{The equations}

We consider a smooth bounded domain $\Omega\subseteq{\bf R}^3$ with $\Gamma=\partial\Omega$ split into $\Gamma_1\cup\Gamma_2$ (with $\Gamma_i$ disjoint subset, $\Gamma_1$ with strictly positive measure).

We assume  that a generalized version of the principle of virtual powers holds, accounting for  internal microforces 
responsible for phase transitions. Thus,
the first principle of thermodynamics reads as follows
\begin{equation}\label{enI}
e_t+\hbox{div }{\bf q }=r+\stress:\strain_t+B_M\chi_{Mt}+B_S\chi_{St}+{\bf B}:\et_t\quad\hbox{in }\Omega
\end{equation}
the right hand side being the power of interior forces and the heat source $r$.
Here, $e$ is the internal energy, ${\bf q}$ the heat flux, $(B_M,B_S)$ and ${\bf B}$ 
internal (microscopic) forces responsible for the phase transformation (i.e. the evolution of internal variables). The heat flux  satisfies boundary condition ($h$ is a known flux through the boundary)
\begin{equation}
{\bf q}\cdot{\bf n}=h\quad\hbox{on }\Gamma.
\end{equation}

Hence, by the principle of virtual powers we get  the quasi-static momentum balance
\begin{equation}
-\hbox{div }\stress={\bf f}\quad\hbox{in }\Omega
\end{equation}
with boundary condition
\begin{align}\label{boundu1}
&\uu={\bf 0}\quad\hbox{on }\Gamma_1\\\label{boundu2}
&\stress{\bf n}={\bf t}\quad\hbox{on }\Gamma_2
\end{align}
${\bf f}$ being a volume force, while ${\bf t}$ is a traction applied on a part of the boundary.

Analogously, the evolution of the phases depends on internal forces which are  
included in the energy balance of the system. 
Thus, we get two balance equations, one for the evolution of the phase proportions (related to $(B_M,B_S)$) 
and one for the evolution of the tensor $\et$ (related to ${\bf B}$), i.e. 
\begin{align}
&(B_M,B_S)=(0,0)\quad\hbox{in }\Omega\\
&{\bf B}={\bf 0}\quad\hbox{in }\Omega.
\end{align}

\subsection{The constitutive relations}

We need to prescribe constitutive relations for the involved physical quantities. 
The internal energy is 
$$
e=\Psi-\theta{\eta}
$$
where the entropy $\eta$ is prescribed by
\begin{align}
&\eta=-\frac{\partial\Psi}{\partial\theta}=c_s\log\theta-h'_A(\theta)(1-\chi_M-\chi_S)\\\no
&-h_M'(\theta)\chi_M-h_S'(\theta)\chi_S-{\bf h}'_d(\theta):\et.
\end{align}
The Cauchy stress tensor is
$$
\stress={\bf S}+\sigma_m{\bf I},
$$
with
$$
\sigma_m=\frac{\partial\Psi}{\partial \hbox{tr }\strain}=\left(\lambda+{\frac 2 3}\mu\right)\hbox{tr }\strain,
$$
and
\begin{equation}\label{Sgrande}
{\bf S}=\frac{\partial\Psi}{\partial{\bf e}}=2\mu({\bf e}-\chi_S\et).
\end{equation}
Hence, we get
\begin{equation}
(B_M,B_S)=(B_M^{nd},B_S^{nd})+(B_M^d,B_S^d)=\frac{\partial\Psi}{\partial(\chi_M,\chi_S)}
+\frac{\partial\phi}{\partial(\chi_{Mt},\chi_{St})}.
\end{equation}
More precisely,
letting 
$$
s(x)=\frac x{|x|}\hbox{ if }x\not=0,\quad s(0)=[-1,1],
$$
and 
$$
H(x)=1\hbox{ if }x>0,\,H(x)=0\hbox{ if }x<0,\,H(0)=[0,1]
$$
there holds
\begin{align}\label{scelgoB1}
&B_M^{nd}=-h_A(\theta)+h_M(\theta)+\frac{\partial {\Psi_{int}}}{\partial\chi_M}+\gamma_M,\\\no
&B_M^d=s(\chi_{Mt}),
\end{align}
and (choosing $\phi_S$ as in  \eqref{specifico})
\begin{align}\label{scelgoB2}
&B_S^{nd}=-h_A(\theta)+h_S(\theta)+\frac{\partial {\Psi_{int}}}{\partial\chi_S}-2\mu(\straind-\chi_S\et):\et+\gamma_S,\\\no 
&B_S^d=\frac{\partial \phi_S}{\partial\chi_{St}}=s(\chi_{St})+d(\chi_S,\stress)H(\chi_{St})
\end{align}
with
\begin{equation}
(\gamma_M,\gamma_S)\in\partial I_K(\chi_M,\chi_S).
\end{equation}
Finally, we consider
\begin{equation}
{\bf B}={\bf B}^{nd}+ {\bf B}^d=\frac{\partial\Psi}{\partial\et}+\frac{\partial\phi}{\partial\et_t}
\end{equation}
where
\begin{align}
&{\bf B}^{nd}=-2\mu\chi_S(\straind-\chi_S\et)+\ttau(\theta)+\gamma\et\\\no
&{\bf B}^d=\chi_S{\bf s}(\et_t)
\end{align}
using the notation
$$
\mathbf s(\et_t)=\frac{\et_t}{\|\et_t\|}\hbox{ if }\et_t\not=0,\quad\mathbf s({\bf 0})=\{{\bf w}:\|{\bf w}\|\leq 1\}.
$$
and letting
$$
\gamma\in\frac 1{\xi_s}\partial I_{\xi_s}(\|\et\|)=\partial I_{\xi_s}(\|\et\|).
$$
As far as the heat flux,  we assume (Fourier law)
\begin{equation}
{\bf q}=-\theta\frac{\partial\Phi}{\partial\nabla\theta}=-\nabla\theta.
\end{equation}

\section{The PDE system}

\subsection{The first principle}

Combining constitutive relations with the balance laws, we get the PDE system we deal with. First let us discuss the energy balance, from which we show that the model is thermodynamically consistent.
The equation governing the evolution of the temperature is recovered from \eqref{enI}. After using the chain rule and by the constitutive relations, we get
\begin{equation}
\theta(\eta_t+\hbox{div }\frac{\bf q}\theta)-r=\frac{\partial\Phi}{\partial(\chi_{Mt},\chi_{St})}\cdot(\chi_{Mt},\chi_{St})+\frac{\partial\Phi}{\partial\et_t}:\et_t+\frac{\partial\Phi}{\partial\nabla\theta}\cdot\nabla\theta\geq0
\end{equation}
from which the second principle of thermodynamics follows, once $\theta>0$ (it is the absolute temperature). Note in particular that we have strongly exploited the fact that $\partial\Phi$ turns out to be a maximal monotone operator with ${\bf 0}\in\partial\Phi({\bf 0})$.
The resulting equation is 
\begin{align}
 &\theta_t(c_s-\theta(h_M''(\theta)\chi_M+h_S''(\theta)\chi_S+h_A''(\theta)(1-\chi_M-\chi_S)+\ttau''(\theta):\et))\\\no
&+\theta h_A'(\theta)(\chi_M+\chi_S)_t-\theta h_M'(\theta)\chi_{Mt}-\theta h_S'(\theta)\chi_{St}-\theta\ttau'(\theta):\et_t-\Delta\theta\\\no
&=|\chi_{Mt}|+|\chi_{St}|+d(\chi_S,\stress)(\chi_{St})^+|\chi_{St}|
+\chi_S\|\et_t\|.
\end{align}
In particular, we have to assume that $h_A,h_M,h_S,{\bf h}_d$ are smooth functions such that
$$
(c_s-\theta(h_M''(\theta)\chi_M+h_S''(\theta)\chi_S+h_A''(\theta)(1-\chi_M-\chi_S)+\ttau''(\theta):\et))\geq C>0.
$$

\subsection{The evolution}

Combining constitutive relations with momentum balance, it follows 
\begin{equation}
-\hbox{div }((\lambda+\frac 2 3\mu)\hbox{tr }\strain{\bf I}+2\mu(\straind-\chi_S\et))={\bf f},
\end{equation}
combined with \eqref{boundu1}, \eqref{boundu2}.
Then, by definition of $B_M$ and $B_S$, the  evolution equations for $(\chi_{Mt},\chi_{St})$ are written as 
\begin{equation}\label{ev1phase}
 s(\chi_{Mt})+(h_M(\theta)-h_A(\theta))+\frac{\partial {\Psi_{int}}}{\partial\chi_M}+\gamma_M=0
\end{equation}
and
\begin{equation}\label{ev2phase}
 s(\chi_{St})+d(\chi_S,\stress)H(\chi_{St})+(h_S(\theta)-h_A(\theta))-2\mu(\straind-\chi_S\et):\et+\gamma_S=0
\end{equation}
where 
$$(\gamma_M,\gamma_S)\in\partial I_K(\chi_M,\chi_S).$$
Note that 
$\partial I_K(\chi_M,\chi_S)=(0,0)$ if $(\chi_M,\chi_S)$ belongs to the interior of $K$, while it is given by the normal cone to the boundary if $(\chi_M,\chi_S)\in\partial K$.

Finally, the evolution equation for $\et$ is given by
\begin{align}\label{ev3phase}
 &\chi_S\mathbf s(\et_t)-2\mu\chi_S(\straind-\chi_S\et)+\ttau(\theta)\\\no
&+\frac{\partial {\Psi_{int}}}{\partial\chi_S}+\gamma\et={\bf 0},\quad\gamma\in\partial I_{\xi_s}(\|\et\|).
\end{align}

\begin{remark}
Note that the coefficient $\chi_S$ of the evolution term $\mathbf s(\et_t)$ ensures that in the absence 
of detwinned martensite there is not dissipative contributions involving  $\et$.
\end{remark}

\subsection{An equivalent formulation}

Let us now introduce ($B_M$ and $B_S$ are defined as in \eqref{scelgoB1} and \eqref{scelgoB2})
\begin{align}
&F_M(B_M^{nd})=|B_M^{nd}|-1\\
&F_S(B_S^{nd})=|B_S^{nd}|-R(B_S^{nd},\chi_S,\stress)
\end{align}
where 
\begin{align}\no
&R(B_S^{nd},\chi_S,\stress)=1\hbox{ if }B_S^{nd}<0\\\no
&\hbox{and }R(B_S^{nd},\chi_S,\stress)=1+d(\chi_S,\stress)\hbox{ if }B_S^{nd}\geq0.
\end{align}
Then, we can rewrite the evolution of the phases \eqref{ev1phase} and \eqref{ev2phase} as follows
\begin{align}
 &\chi_{Mt}=\zeta_M\frac{B_M^{nd}}{|B_M^{nd}|},\\
&\chi_{St}=\zeta_S\frac{B_S^{nd}}{|B_S^{nd}|},\\
&\zeta_i F_i=0\quad i=M,S,\\
&F_i(B_i^{nd})\leq 0,\quad i=M,S.
\end{align}
Note that ${F}_i$ play the role of yield functions (see, e.g.,  \cite{lub}).

Analogously we may introduce
$$
{F_d}=\|{\bf B}^{nd}\|-\chi_S
$$
letting
$$
\et_t=\zeta_d\frac{{\bf B}^{nd}}{\|{\bf B}^{nd}\|},
$$
with
$$
\zeta_d{F_d}=0,\qquad F_d\leq 0.
$$

\section{Some examples}

In the following we explore the model performances limiting the discussion only to the case of a proportional loading state, i.e., neglecting the reorientation process. 
Accordingly, to simplify the discussion, we set
$\ttau(\theta)={\bf 0}$, $\et$ in the same direction of ${\bf S}$ (and $\straind$), assuming also $\et_t=  {\bf 0}$.
Under these simplifying positions, we may set (see \eqref{Sgrande})
\begin{equation}
\frac{{\bf B}^{nd}}{||{\bf B}^{nd} ||} =
\frac{\straind}{||\straind ||}  =
\frac{\et}{|| \et ||} ,
\end{equation}
and 
\begin{equation}
|| {\bf S} || = 2\mu ( || {\bf e} || -\chi_S).
\end{equation}
Moreover, we distinguish between two different possible situations, one in which we consider only evolution of the stress-induced martensite and one in which we consider only evolution of the temperature-induced martensite, as discussed in the following.
For both problems we start from a material completely in austenite (i.e., 
$\chi_S=\chi_M=0$).

\subsection{Case 1: temperature-induced effect}

For this problem we assume to start from 
$\sigma=0$ and to vary only the temperature. Accordingly,
only a variation of $\chi_M$ can be produced.

The problem is governed by the following set of equations (see \eqref{scelgoB1}):
\begin{equation}
\left\{
\begin{aligned}
& B_M^{nd}=-h_A(\theta)+h_M(\theta)+\frac{\partial {\Psi_{int}}}{\partial\chi_M}+\gamma_M,
\\
&F_M(B_M^{nd})=|B_M^{nd}|-1
\\
&\chi_{Mt}=\zeta_M\frac{B_M^{nd}}{|B_M^{nd}|},
\\
&\zeta_M F_M=0 \quad F_M(B^{nd}_M) \leq 0 .
\end{aligned}
\right.
\end{equation}
We assume to first properly cool and then heat the material (see Figure 
\ref{chiM_sigma_T}).

In Figures 
\ref{chiM_B_T}-\ref{chiM_gammaM_T}
we report the evolution
of the thermodynamic force $B_M^{nd}$ versus the temperature $\theta$,
of the temperature-induced martensite $\chi_S$ versus the temperature $\theta$,
of the quantity $\gamma_M$ versus the the temperature $\theta$.

It can be observed that during cooling the model is able to reproduce a process in which  the multi-variant martensite is produced and then during heating a process in which  the multi-variant martensite is progressively extinguished. The forward and reverse phase transformations are perfectly symmetric.

\subsection{Case 2: stress-induced effect}

For this problem we assume to start from 
$\sigma=0$  and to vary only the stress. Accordingly,
only a variation of $\chi_S$ can be produced.

The problem is governed by the following set of equations (\eqref{scelgoB2}):
\begin{equation}
\left\{
\begin{aligned}
& B_S^{nd}=-h_A(\theta)+h_S(\theta)+\frac{\partial {\Psi_{int}}}{\partial\chi_S}-
2\mu(|| \straind||-\chi_S)+\gamma_S,
\\
& F_S(B_S^{nd})=|B_S^{nd}|-R(B_S^{nd},\chi_S,\stress)
\\
& \chi_{St}=\zeta_S\frac{B_S^{nd}}{|B_S^{nd}|},
\\
&\zeta_S F_S=0 \quad F_S(B^{nd}_S) \leq 0 .
\end{aligned}
\right.
\end{equation}
where (letting \eqref{specifico} holds)
\begin{equation}
\begin{aligned}
& R(B_S^{nd},\chi_S,\stress)=1\hbox{ if }B_S^{nd}<0\hbox{ and }
\\
& R(B_S^{nd},\chi_S,\stress)=1+d(\chi_S,\stress)\hbox{ if }B_S^{nd}\geq0.
\end{aligned}
\end{equation}

We assume to first properly load and then unload the material (see Figure 
\ref{chiS_sigma_T}).

In Figures 
\ref{chiS_B_sigma}-\ref{chiS_sigma_eps}
we report the evolution
of the thermodynamic force $B_S^{nd}$ versus the applied stress $	\sigma$,
of the stress-induced martensite $\chi_S$ versus the applied stress $	\sigma$,
of the quantity $\gamma_S$ versus the applied stress $\sigma$,
of  the applied stress $\sigma$ versus the strain $\epsilon$.

It can be observed that during loading the model is able to reproduce a process in which  the single-variant martensite is produced and then during unloading a process in which  the single-variant martensite is progressively extinguished. The forward and reverse phase transformation are unsymmetric.

ù

\end{document}